\documentclass{article}

\usepackage[T2C]{fontenc}
\usepackage[tbtags]{amsmath}
\usepackage{amsfonts,amssymb,mathrsfs,amscd}
\usepackage{hyperref}

\usepackage{theorem}

\newtheorem{theorem}{Theorem}[section]

\newtheorem{proposition}[theorem]{Proposition}

\newtheorem{proof}{Proof}
{\theorembodyfont{\rmfamily}
\newtheorem{definition}{Definition}
\newtheorem{remark}[theorem]{Remark}
\newtheorem{example}[theorem]{Example}
}

\allowdisplaybreaks

\def\Res{\operatornamewithlimits{Res}}

\newcounter{aa}

\let\mathcal=\mathscr

\begin{document}



\author{M.~Kazarian, S.~Lando}

\date{}

\title{Combinatorial solutions to integrable hierarchies}

\maketitle
\markright{Combinatorial solutions to integrable hierarchies}

\footnotetext[0]{The work was partly supported by the RFBR
(grant 
~13-01-00383-a).}

\begin{abstract}
We give a review of modern approaches to constructing formal solutions to
integrable hierarchies of mathematical physics, whose coefficients
are answers to various enumerative problems. The relationship between these
approaches and combinatorics of symmetric groups and their representations is explained.
Applications of the results to constructing efficient computations in
problems related to models of quantum field theories are given.

\end{abstract}


\setcounter{tocdepth}{2} \tableofcontents

\section*{Introduction}

One of the models of two-dimensional quantum gravity suggests that
a random metric surface should be discretized, that is, represented
as a result of gluing flat polygons. In particular, one can think
about gluings of triangles (triangulations) or squares (quadrangulations).
A study of such a model requires analyzing asymptotics
of various combinatorial enumerative sequences related to gluings
(e.g., the asymptotics of the number of triangulations of a surface
of a given genus as the number of triangles tends to infinity).
In spite of serious efforts both by mathematicians and physicists
during the last decades, application of the methods developed
to problems of this kind is still not automatic and requires
nontrivial skills. Simultaneously, a large class of enumerative problems
was discovered, the answers to which are solutions to integrable hierarchies of
mathematical physics.

Appearance of integrable hierarchies (namely, the Korteweg--de~Vries
hierarchy) in this kind of problems is known for more than thirty years;
however, then existing constructions required a nontrivial study of
singularities of generating functions in question and the so-called
double-scaling limit analysis. The first natural example is the generating
function for the double Hurwitz numbers. R.~Pandharipande~\cite{25}
conjectured in~2000 that this function is a solution to the Toda lattice
hierarchy, and A.~Okounkov~\cite{22} proved the conjecture the same year
(see Sec.~\ref{ssec2.4}). In its own turn, the interest to the study of Hurwitz
numbers arose from that to the geometry of moduli spaces of complex curves.
This geometry serves as an alternative module of two-dimensional quantum gravity,
and its relationship to the Korteweg--de~Vries hierarchy was conjectured
by E.~Witten~\cite{30} basing on the double scaling limit.
Okounkov's result allowed the present authors to give an algebro-geometric
proof of Witten's conjecture~\cite{16}.

After Okounkov's work, natural combinatorial solutions to integrable hierarchies
began to multiply. It happened that nonlinear partial differential equations
can serve as efficient computational tools in enumerating various combinatorial
objects, which, in turn, allows one to compute various asymptotics of the
corresponding sequences that were not known earlier.

The present paper is divided into two parts. In the first part, we present
a general description of the Kadomtsev--Petviashvili integrable hierarchy,
which can serve as a standard sample hierarchy of this kind. We describe a
general construction of solutions to this hierarchy, due to Sato.

The second part deals with a certain special infinite-dimensional family of solutions
to the Kadomtsev--Petviashvili hierarchy, and a series of examples of enumerative
generating functions belonging to this family is given. Then we show how these solutions
are related to enumeration of polygonal gluings. We only mention relationship
of the solutions thus constructed with the geometry of moduli spaces of complex curves.

The authors are grateful to participants of the ``Characteristic classes
and intersection theory'' seminar at the Math department of the Higher School of
Economics for numerous fruitful discussions.

\section{Kadomtsev--Petviashvili hierarchy}
\label{sec1}

The Kadomtsev--Petviashvili hierarchy (below, KP hierarchy, for brevity) is a totally integrable system
of partial differential equations for functions depending on infinitely many variables.
The physical origins of this hierarchy are well described, for example, in~\cite{19},
and we do not touch on them. We are interested only in certain functions satisfying all
the equations of this hierarchy, and we will use a description of all its solutions
presented below.

\subsection{Equations of the KP hierarchy}
\label{ssec1.1}
Mathematical physicists use different notation for variables, depending
on the chosen normalization; we will use the variables
$p_1,p_2,\dots$, such that the lowest KP hierarchy equation has the form
\begin{equation}
\label{eq1}
\frac{\partial^2 F}{\partial p_2^2}=
\frac{\partial^2 F}{\partial p_1\,\partial p_3}-\frac{1}{2}
\biggl(\frac{\partial^2 F}{\partial p_1^2}\biggr)^2-
\frac{1}{12}\,\frac{\partial^4 F}{\partial p_1^4}\,,
\end{equation}
or
$$
F_{2^2}=F_{1^13^1}-\frac{1}{2}(F_{1^2})^2-\frac{1}{12}F_{1^4}.
$$
(Here the partition in the lower index of a function consists of the indices
of the variables~$p$ with respect to which the partial derivative is taken.
We make use of the multiplicative representation of partitions:
thus, $2^23^1$ denotes the partition of~$7$ having two parts equal to~$2$
and the third part equal to~$3$.) The equation contains derivatives with respect to
the first three variables only, and the sum of the indices in each monomial is~$4$.
The indices in the next equations contain partitions of $5,6,\dots$:
\begin{align*}
F_{2^13^1}&=-F_{1^2}F_{1^12^1}+F_{1^14^1}-\frac{1}{6}F_{1^32^1},
\\
F_{2^14^1}&=-\frac{1}{2}F^2_{1^12^1}-F_{1^2}F_{1^13^1}+
\frac{1}{8}F^2_{1^3}+\frac{1}{12}F_{1^2}F_{1^4}
-\frac{1}{4}F_{1^33^1}+\frac{1}{120}F_{1^6},
\\
F_{3^2}&=\frac{1}{3}F^3_{1^2}-F^2_{1^12^1}-F_{1^2}F_{1^13^1}+
F_{1^15^1}+\frac{1}{4}F^2_{1^3}+\frac{1}{3}F_{1^2}F_{1^4}-
\frac{1}{3}F_{1^33^1}+\frac{1}{45}F_{1^6}.
\end{align*}
Note that the left-hand sides of the equations are second derivatives of
the unknown function with respect to variables having numbers greater than~$1$,
while all the derivatives but, possibly, one on the right-hand side are with respect to
the variable~$p_1$.

\subsection{$\tau$-functions of the KP hierarchy}
\label{ssec1.2}
Exponents of solutions to the KP hierarchy are called the $\tau$-functions
of the hierarchy. Let us give an example of a $\tau$-function
(whose logarithm is, of course, a solution to KP). In order to define
its coefficients, we will require the following definition.

\begin{definition}
\label{def1}
Let $\mu$ be a partition, $\mu\vdash |\mu|$.
The \textit{simple Hurwitz number~$h^\circ_{m;\mu}$} is defined by the equation
$$
h^\circ_{m;\mu}=\frac{1}{|\mu|!}|\{(\eta_1,\dots,\eta_m),
\eta_i\in C_2(S_{|\mu|})\colon \eta_m\circ\dots\circ
\eta_1\in C_\mu(S_{|\mu|})\}|.
$$
Here $S_{|\mu|}$ is the group of permutations of~$|\mu|$ elements,
$C_2(S_{|\mu|})$ denotes the set of all transpositions in~$S_{|\mu|}$,
and $C_\mu(S_{|\mu|})$ is the set of all permutations of cyclic type~$\mu\vdash |\mu|$ in~$S_{|\mu|}$.
(The \textit{cyclic type} of a permutation is the set of lengths of the cycles
in its decomposition in a product of independent cycles;
it forms a partition of the number of permuted elements). In particular,
$$
C_2(S_{|\mu|})=C_{1^{|\mu|-2}2^1}(S_{|\mu|}).
$$
In other words, simple Hurwitz numbers enumerate decompositions of
a given permutation into a product of a given number of transpositions,
in terms of the cyclic type of the permutation.

The \textit{connected simple Hurwitz number~$h_{m;\mu}$} is defined similarly,
but only those sequences of~$m$ permutations are taken into account such that
the subgroup $\langle\eta_1,\dots,\eta_m\rangle\subseteq S_{|\mu|}$
generated by them acts transitively on the set $\{1,\dots,|\mu|\}$.
\end{definition}

Collect the simple Hurwitz numbers into two generating functions:
\begin{align}
\label{eq2}
H^\circ(u;p_1,p_2,\dots)&=\sum_{m=0}^\infty\,\sum_\mu
h^\circ_{m;\mu}p_{\mu_1}p_{\mu_2}\cdots\frac{u^m}{m!}\,,
\\
\label{eq3}
H(u;p_1,p_2,\dots)&=\sum_{m=1}^\infty\,\sum_{\mu\ne\varnothing}
h_{m;\mu}p_{\mu_1}p_{\mu_2}\cdots\frac{u^m}{m!}\,,
\end{align}
where in the second case~$\mu$ runs over all partitions of all natural numbers,
while in the first case the empty partition also is added (which corresponds
to the summand~$1$). These generating functions depend on infinitely many
variables and are formal: we impose no restrictions on their convergence.
It is clear from the definition that the coefficients of both~$H$ and~$H^\circ$
are rational. The first terms of the power series expansions are
\begin{align*}
H^\circ(u;p_1,p_2,\dots)&=1+p_1+\frac{p_2 u}{2}+\frac{p_1^2}{2}+
\frac{p_3 u^2}{2}+\frac{1}{2} p_2 p_1 u+\frac{p_1^3}{6}+
\frac{p_2 u^3}{12}
\\
&\qquad+\frac{2 p_4 u^3}{3}+\frac{1}{4} p_1^2 u^2+
\frac{1}{2} p_3 p_1 u^2+\frac{1}{8} p_2^2 u^2+\frac{1}{4} p_2 p_1^2 u+
\frac{p_1^4}{24}+\cdots,
\\
H(u;p_1,p_2,\dots)&=p_1+\frac{p_2 u}{2}+\biggl(\frac{p_1^2}{4}+
\frac{p_3}{2}\biggr) u^2+\biggl(\frac{2 p_1 p_2}{3}+\frac{p_2}{12}+
\frac{2 p_4}{3}\biggr)u^3
\\*
&\qquad+\biggl(\frac{p_1^3}{6}+\frac{p_1^2}{48}+\frac{9 p_3 p_1}{8}+
\frac{p_2^2}{2}+\frac{3 p_3}{8}+\frac{25p_5}{24}\biggr)u^4+\cdots\,.
\end{align*}

A general combinatorial statement about enumeration of connected and
disconnected objects yields
$$
H^\circ=\exp\{H\}.
$$
As a result, one can reformulate properties of simple
Hurwitz numbers into those of connected simple ones and vice versa.

\begin{theorem}[\cite{22},~\cite{16}]
\label{th1.1}
The generating function~$H^\circ(u;p_1,p_2,\dots)$ for simple Hurwitz numbers
is a one-parametric family of $\tau$-functions for the KP hierarchy, while
its logarithm~$H(u;p_1,p_2,\dots)$ is a one-parametric family of solutions to this hierarchy.
\end{theorem}

The Sato construction~\cite{27} below interprets the
space of solutions of the KP hierarchy as the Grassmanian of semi-infinite
planes in an infinite-dimensional vector space. This construction allows
one to produce all the solutions of the KP equations.
In order to solve the inverse problem, that is, to prove that a given function
is a solution to KP, it suffices to specify the corresponding semi-infinite plane.
This means, in particular, that in order to identify a function as a solution
to KP hierarchy, there is no need to know what are the equations of the hierarchy.
Below, we will, among others, specify the semi-infinite plane corresponding to the Hurwitz
function~$H$.

\subsection{Projective embeddings of Grassmanians and Pl\"ucker coordinates}
\label{ssec1.3}
Consider the Grassmanian~$G(2,V)$ of two-dimensional planes
in the $4$-dimensional vector space $V\equiv\mathbb{C}^4$.
Any two-dimensional plane in~$V$ can be represented as a wedge product
$\beta_1\wedge \beta_2$ of any pair of linearly independent vectors~$\beta_1$,~$\beta_2$ in it.
This wedge product is well-defined up to a nonzero factor; it determines
the two-dimensional plane uniquely. Hence, this wedge product defines an
embedding of the Grassmanian~$G(2,V)$ into the projectivization
of the wedge square of~$V$:
$G(2,V)\hookrightarrow P\Lambda^2V$.
This construction admits an immediate generalization to embeddings of an arbitrary
Grassmanian~$G(k,V)$ of
$k$-dimensional planes in an~$n$-dimensional vector space~$V$
to the projectivization~$P\Lambda^kV$.

The \textit{Pl\"ucker equations} are the equations of the image of
this embedding. For $\dim V=n$, the Grassmanian~$G(k,V)$ has dimension~$k(n-k)$,
while the dimension of~$P\Lambda^k V$ is
$\begin{pmatrix} n\\ k\end{pmatrix}-1$,
which, as a rule, is greater than~$k(n-k)$;
hence, the image of the Pl\"ucker embedding generally does not coincide
with the projectivized  wedge product~$P\Lambda^kV$.
For example, the image of the embedding of~$G(2,4)$ into~$P\Lambda^2\mathbb{C}^4$
is a hypersurface in a $5$-dimensional projective space.

In order to find the equation of this hypersurface, let us pick a
basis~$e_1$, $e_2$, $e_3$,~$e_4$ in~$V$. This basis generates the basis
$\beta_{ij}=e_i\wedge e_j$, $1\leqslant i<j\leqslant4$,
in~$\Lambda^2 V$ and the natural coordinate system~$y_{ij}$
in the latter vector space. The image of the embedding of
the Grassmanian consists of decomposable vectors. By the definition of the
wedge product, the image of the plane spanned by a pair of vectors
$(a_1,a_2,a_3,a_4)$, $(b_1,b_2,b_3,b_4)$, has projective coordinates
$$
y_{ij}=\begin{vmatrix} a_i & b_i
\\
a_j & b_j\end{vmatrix}=a_ib_j-a_jb_i,
$$
whence
$$
y_{12}y_{34}-y_{13}y_{24}+y_{14}y_{23}=0.
$$
This is exactly the Pl\"ucker equation of the embedding.

For arbitrary values of~$n$ and~$k$, the Pl\"ucker equations remain quadratic.
In other words, the ideal in the ring of polynomials consisting of
polynomials vanishing on the image of the Pl\"ucker embedding
is generated by quadratic polynomials (see, e.g.,~\cite{9}).

\subsection{The vector space of Laurent series and its semi-infinite wedge power}
\label{ssec1.4}
Take for the space~$V$ the infinite dimensional vector space of Laurent series in one variable.
Elements of this space are series of the form
$$
c_{-k}z^{-k}+c_{-k+1}z^{-k+1}+\cdots\,.
$$
The monomials~$z^k$, $k=\dots,-2,-1,0,1,2,\dots$,
form the \textit{standard} basis in~$V$. By definition, the
\textit{semi-infinite wedge product $\Lambda^{\infty/2}V$}
is the vector space freely spanned by the vectors
$$
v_\mu=z^{m_1}\wedge z^{m_2}\wedge z^{m_3}\wedge\cdots,\qquad
m_1<m_2<m_3<\cdots,\quad
m_i=\mu_i-i,
$$
where $\mu$ is a partition, $\mu=(\mu_1,\mu_2,\mu_3,\dots)$,
$\mu_1\geqslant\mu_2\geqslant\mu_3\geqslant\cdots$,
with all but finitely many parts equal to~$0$.
In particular, $m_i=-i$ for all~$i$ sufficiently large.

The \textit{vacuum vector}
$$
v_{\varnothing}=z^{-1}\wedge z^{-2}\wedge z^{-3}\wedge\cdots
$$
corresponds to the empty partition. Similarly,
\begin{align*}
v_{1^1}&=z^{0}\wedge z^{-2}\wedge z^{-3}\wedge\cdots,
\\
v_{2^1}&=z^{1}\wedge z^{-2}\wedge z^{-3}\wedge\cdots,
\\
v_{1^2}&=z^{0}\wedge z^{-1}\wedge z^{-3}\wedge\cdots\quad
\text{and so on.}
\end{align*}

\subsection{Boson--fermion correspondence}
\label{ssec1.5}
Enumeration of basic vectors of the semi-infinite wedge product~$\Lambda^{\infty/2}V$
(that is, of the space of \textit{fermions}) by partitions of nonnegative integers
establishes its natural vector space isomorphism with the space of \textit{bosons},
the space of power series in infinitely many variables $p_1,p_2,\dots$.
This isomorphism takes a basic vector~$v_\mu$ to the Schur polynomial
$s_\mu=s_\mu(p_1,p_2,\dots)$. The latter is a quasihomogeneous polynomial
of degree~$|\mu|$ in the variables~$p_i$
(by definition, the degree of~$p_i$ is~$i$). The Schur polynomials form an additive
basis of the vector space of power series in the variables~$p_i$.

The \textit{Schur polynomial~$s_k$} corresponding to the one-part partition
$k^1\vdash k$ is defined by means of the power series expansion
\begin{align*}
&s_0+s_1z+s_2z^2+s_3z^3+s_4z^4+\cdots=\exp\biggl\{p_1z+p_2\frac{z^2}{2}+
p_3\frac{z^3}{3}+\cdots\biggr\}
\\
&\qquad=1+p_1z+\frac{1}{2}(p_1^2+p_2)z^2+
\frac{1}{6}(p_1^3+3p_1p_2+2p_3)z^3+\cdots.
\end{align*}
For an arbitrary partition $\kappa=(\kappa_1,\kappa_2,\kappa_3,\dots)$,
$\kappa_1\geqslant\kappa_2\geqslant\kappa_3\geqslant\cdots$,
the Schur polynomial~$s_\kappa$ is given by the determinant
\begin{equation}
\label{eq4}
s_\kappa=\det\|s_{\kappa_j-j+i}\|.
\end{equation}
The indices~$i$,~$j$ here run over the set $\{1,2,\dots,n\}$,
for~$n$ large enough; since for~$i$ large enough we have
$\kappa_i=0$, this determinant, whence the polynomial~$s_\kappa$,
is independent of~$n$. Thus,
\begin{align*}
s_{1^3}&=\begin{vmatrix}
s_1&s_2&s_3&s_4&\hdots\\
s_0&s_1&s_2&s_3&\hdots\\
0&s_0&s_1&s_2&\hdots\\
0&0&0&s_0&\hdots\\
\hdotsfor{5}
\end{vmatrix}
\\
&=\begin{vmatrix}
p_1&\dfrac{1}{2}(p_1^2+p_2)&\dfrac{1}{6}(p_1^3+3p_1p_2+2p_3)&\hdots\\
1&p_1&\dfrac{1}{2}(p_1^2+p_2)&\hdots\\
0&1&p_1&\dfrac{1}{2}(p_1^2+p_2)\\
0&0&0&1\\
\hdotsfor{4}
\end{vmatrix}
\\
&=\begin{vmatrix}
p_1&\dfrac{1}{2}(p_1^2+p_2)&\dfrac{1}{6}(p_1^3+3p_1p_2+2p_3)\\
1&p_1&\dfrac{1}{2}(p_1^2+p_2)\\
0&1&p_1
\end{vmatrix}.
\end{align*}

Here are the first few Schur polynomials:
\begin{gather*}
s_0=1,\quad
s_{1^1}=p_1,\quad
s_{2^1}=\frac{1}{2}(p_1^2+p_2),\quad
s_{3^1}=\frac{1}{6}(p_1^3+3p_1p_2+2p_3),
\\
s_{1^2}=\frac{1}{2}(p_1^2-p_2), \quad
s_{1^12^1}=\frac{1}{3}(p_1^3-p_3),\quad
s_{1^3}=\frac{1}{6}(p_1^3-3p_1p_2+2p_3).
\end{gather*}

The Schur polynomials play a key role in the study of representations
of symmetric groups, hence they appear in enumerative problems
for product decompositions of permutations not by chance.

There is a more conceptual description of the boson--fermion correspondence.
Namely, for $i=1,2,3,\dots$, denote by~$\widehat{z^{-i}}$
the shift operator acting in~$V$ and taking~$z^m$ to~$z^{m-i}$, $m\in\mathbb{Z}$.
This action extends to the space of fermions by the Leibnitz rule:
$$
\widehat{z^{-i}}\colon z^{m_1}\wedge z^{m_2}\wedge\cdots\mapsto
z^{m_1-i}\wedge z^{m_2}\wedge\cdots+z^{m_1}\wedge z^{m_2-i}
\wedge\cdots+\cdots\,.
$$
Let~$\langle v\rangle$ denote the coefficient of the vacuum vector
in the decomposition of a vector~$v$ in the basis~$v_\mu$.
In this notation, the boson--fermion correspondence takes a vector
$v\in\Lambda^{\infty/2}V$ to the following formal series in the variables~$p_k$:
$$
v\mapsto \biggl\langle \exp\biggl\{\,\sum_{i=1}^\infty\frac{p_i}{i}
\widehat{z^{-i}}\biggr\}v\biggr\rangle.
$$
This formula can be interpreted as follows. The boson--fermion correspondence
takes the operator~$\widehat{z^{-i}}$ acting in the space of fermions
to the operator $i\partial/\partial p_i$ acting in the space of bosons.
Thus, the formula above is nothing but the Taylor formula
$$
f(p)=\exp\biggl\{\,\sum_{i=0}^\infty p_i\frac{\partial}{\partial q_i}\biggr\}
f(q)\big|_{q=0}.
$$
Note also that for positive~$i$ the shift operator~$\widehat{z^{i}}$
acts in the space of bosons as multiplication by~$p_i$.
Hence, in the opposite direction the boson--fermion correspondence
takes a polynomial (or, more generally, a series) $f(p_1,p_2,\dots)$
to the fermion~$f(\widehat{z^{1}},\widehat{z^{2}},\dots)v_\varnothing$.

\subsection{Semi-infinite Grassmanian and KP equations}
\label{ssec1.6}
The \textit{semi-infinite Grassmanian~$G(\infty/2,V)$}
consists of decomposable elements in~$P\Lambda^{\infty/2}V$,
i.e., of vectors of the form
$$
\beta_1(z)\wedge \beta_2(z)\wedge \beta_3(z)\wedge\cdots,
$$
where $\beta_i$ are Laurent series in the variable~$z$,
and for~$i$ large enough the leading term of the series expansion of~$\beta_i$
is~$z^{-i}$:
$$
\beta_i(z)=z^{-i}+c_{i1}z^{-i+1}+c_{i2}z^{-i+2}+\cdots\,.
$$
Such an infinite wedge product considered up to a nonzero factor
corresponds to the semi-infinite plane spanned by the vectors
$\beta_1,\beta_2,\dots$\,.

\begin{definition}
\label{def2}
The \textit{Hirota equations} are the Pl\"ucker equations for
the embedding of the semi-infinite Grassmanian into the projectivized
wedge product $P\Lambda^{\infty/2}V$. Solutions of the Hirota equations
(that is, semi-infinite planes) are called the
\textit{$\tau$-functions} of the KP hierarchy.
\end{definition}

Here is a way to represent the Hirota equations: the following identities
must be satisfied identically in~$p$ and~$q$:
\begin{align*}
&\Res_{z=0}\exp\,\biggl\{2\sum_{i=0}^\infty
\frac{q_i}{i}z^{-i}\biggr\}
\\
&\qquad\times\tau(p_1-q_1+z,p_2-q_2+z^2,\dots)
\tau(p_1+q_1-z,p_2+q_2-z^2,\dots)\frac{dz}{z^2}=0;
\end{align*}
here $\Res_{z=0}$ denotes the operation of taking the residue
of a~$1$-form in the variable~$z$ at $z=0$.

This identity can be interpreted as follows. Take the formal series expansion
of the left-hand side in the $q$-variables. Then the coefficient of each
monomial is a finite quadratic polynomial combination of partial derivatives
of the function $\tau=\tau(p_1,p_2,\dots)$. For example, the equation
corresponding to the monomial~$q_3$ has the form
$$
\tau_2^2-\tau_1\tau_3-\frac{1}{4}\tau_{1^2}^2+\tau_1\tau_3-\tau\,\tau_{2^2}+
\frac{1}{3}\tau_1\tau_{1^3}-\frac{1}{12}\tau\,\tau_{1^4}=0.
$$
Each of these quadratic differential equations can be
treated as a system of quadratic algebraic relations
in the Taylor coefficients of~$\tau$, which are exactly the Pl\"ucker equations.
On the other hand, the substitution $\tau=e^F$ transforms the equations
for the function~$\tau$ into equations for the function~$F$.
For example, the equation above is equivalent to Eq.~\eqref{eq1}.

\begin{definition}
\label{def3}
The form the Hirota equations acquire for logarithms of
$\tau$-functions is called the
\textit{Kadomtsev--Petviashvili hierarchy}.
\end{definition}

In other words, each solution to the KP hierarchy can be obtained as the result of
the following procedure:

-- take the wedge product $\beta_1(z)\wedge\beta_2(z)\wedge\cdots$
corresponding to some semi-infinite plane in~$V$;

-- taking the Laurent series expansion of the functions~$\beta_i$, represent
the corresponding point in the semi-infinite Grassmanian as a linear combination
(which contains, generally speaking,
infinitely many nonzero summands)
of the basic vectors~$v_\kappa$  and multiply it by a constant in order
to make the coefficient of the vacuum vector~$v_\varnothing$ equal to~$1$;

-- replace each vector~$v_\kappa$ in the resulting linear combination
by the corresponding Schur polynomial~$s_\kappa(p_1,p_2,\dots)$
thus obtaining a power series in infinitely many variables
$p_1,p_2,\dots$;

-- take the logarithm of the resulting series.

Thus we have a universal tool to produce all the solutions to the KP hierarchy.
Checking whether a given function is a solution to the KP hierarchy usually is
more complicated. We construct the semi-infinite plane corresponding to the
generating function~$H^\circ$ for simple Hurwitz numbers in Sec.~\ref{ssec2.2}.

\section{
Combinatorial representatives of the Orlov--Shcherbin family of solutions
to the KP hierarchy}
\label{sec2}

The generating function for the Hurwitz numbers belongs to a vast family
of solutions to the KP hierarchy. This family contains a lot of other
combinatorial solutions whose coefficients enumerate various objects related
to combinatorics of symmetric groups.
In this section we describe this family and indicate several specializations of its
parameters giving answers to some enumerative problems.

\subsection{Description of the family}
\label{ssec2.1}
To each partition $\mu=(\mu_1,\dots,\mu_\ell)$ of~$n$,
a Young diagram consisting of~$n$ cells is associated;
we denote the number~$\ell$ of rows in this diagram, that is, the number of
parts in the partition~$\mu$ by~$\ell(\mu)$. The \textit{content} of a cell~$w$
of a Young diagram belonging to the $i$~th row
and the $j$~th column is the quantity
$$
c(w)=j-i, \qquad
1\leqslant i\leqslant\ell(\mu),\quad
1\leqslant j\leqslant\mu_i.
$$
The \textit{content~$c(\mu)$ of a diagram~$\mu$}
is the unordered tuple of the contents of its cells. Set also
$$
y_\mu=\prod_{w\in\mu}y_{c(w)}.
$$
For example,
$$
c((5,3,3,2))=(-3,-2,-2,-1,-1,0,0,0,1,1,2,3,4)
$$
and
$$
y_{(5,3,3,2)}=y_{-3}y_{-2}^2y_{-1}^2y_0^3y_1^2y_2y_3y_4;
$$
and the degree of this monomial is the number of cells in the diagram
(this number equals~$13$ in the example, see Fig.~\ref{fig1}).

\begin{figure}[h]
\centerline{\begin{picture}(100,90)
\thicklines
\put(0,0){\line(1,0){40}}
\multiput(0,20)(0,20){2}{\line(1,0){60}}
\multiput(0,60)(0,20){2}{\line(1,0){100}}
\multiput(0,0)(20,0){3}{\line(0,1){80}}
\put(60,20){\line(0,1){60}}
\multiput(80,60)(20,0){2}{\line(0,1){20}}
\multiput(7,67)(20,-20){3}{$0$}
\multiput(27,67)(20,-20){2}{$1$}
\put(47,67){$2$}
\put(67,67){$3$}
\put(87,67){$4$}
\multiput(2,47)(20,-20){2}{$-1$}
\multiput(2,27)(20,-20){2}{$-2$}
\put(2,7){$-3$}
\end{picture}}
\caption{A Young diagram, with the content written inside each cell}
\label{fig1}
\end{figure}

We present below the theorem of A.~Orlov and D.~Shcherbin~\cite{24}
in the form due to I.~Goulden and D.~Jackson~\cite{8}.

\begin{theorem}
\label{th2.1}
The generating function
$$
\sum_\mu y_\mu\frac{\dim_\mu}{|\mu|!}s_\mu(p_1,p_2,\dots),
$$
where $\dim_\mu$ is the dimension of the irreducible representation of
the symmetric group~$S_{|\mu|}$ corresponding to a partition~$\mu$,
is a family of $\tau$-functions for the KP hierarchy.
\end{theorem}

Note that the factor $\dim_\mu/|\mu|!$ can be written in the form
$$
\frac{\dim_\mu}{|\mu|!}=s_\mu(1,0,0,\dots).
$$

In particular, substituting for the variables
$\dots,y_{-2},y_{-1},y_0,y_1,y_2,\dots$ any specific values
we obtain concrete $\tau$-functions, and their logarithms
yield solutions to the KP hierarchy.

\begin{proof}
By setting $y_c\equiv1$, $c=\dots,-2,-1,0,1,2,\dots$,
we make the $\tau$-function above into
$$
\sum_\mu \frac{\dim_\mu}{|\mu|!}s_\mu(p)=e^{p_1},
$$
which is a specialization of the following well-known identity for the Schur polynomials:
$$
\exp\biggl\{\,\sum_{k=1}^\infty\frac{p_kq_k}{k}\biggr\}=
\sum_\mu s_\mu(p)s_\mu(q)
$$
(see, e.g.,~\cite{26}). This $\tau$-function corresponds, in the fermion
representation, to the infinite wedge product
$\beta_1\wedge\beta_2\wedge\beta_3\wedge\dotsb$, where
\begin{equation}
\label{eq5}
\beta_k=e^zz^{-k}=\sum_{i=0}^\infty\frac{z^{i-k}}{i!}\,.
\end{equation}
Note that its logarithm~$p_1$ is a solution to the KP hierarchy because of the
obvious reason: all the derivatives in the equations of the KP hierarchy
are of order at least~$2$.

For general values of the parameters~$y_c$, the latter fermionic
representation must be modified as follows. Set
$$
u_k=\begin{cases}
\displaystyle\prod_{i=1}^ky_i,&k>0,
\\
\biggl(\,\displaystyle\prod_{i=k+1}^0 y_i\biggr)^{-1},&k\leqslant0.
\end{cases}
$$
In other words, we have chosen~$u_k$ in such a way that, for all integer $k\leqslant\ell$
the equation $\dfrac{u_k}{u_\ell}=\displaystyle\prod_{i=k+1}^\ell y_i$ holds.
In this notation, the factors in the wedge product corresponding to the
$\tau$-function in the theorem are
$$
\beta_k=\frac{1}{u_{-k}}\sum_{i=0}^\infty u_{i-k}\frac{z^{i-k}}{i!}\,.
$$
Indeed, under such a modification of the wedge product, the coefficient
of $v_\mu=z^{\mu_1-1}\wedge z^{\mu_2-2}\wedge\cdots$
in its expansion in basic vectors will be multiplied by the monomial
$\displaystyle\prod_{k=1}^\infty\dfrac{u_{\mu_k-k}}{u_{-k}}=y_\mu$, as required.
The theorem is proved.
\end{proof}

\begin{remark}
\label{rem2.2}
Acting on a semi-infinite plane in the space of Laurent series~$V$
by a linear operator we generally obtain a different semi-infinite plane,
which, in its own turn, produces a different $\tau$-function of the
KP hierarchy. In particular, the $\tau$-functions in the Orlov--Shcherbin family
are obtained from the $\tau$-function~$e^{p_1}$ as a result of the action on the plane~\eqref{eq5}
of the linear operator whose matrix in the standard basis is
$$
\begin{pmatrix}
\ddots&\hdots&\hdots&\hdots&\hdots&\hdots&\hdots&\hdots
\\
\hdots&\dfrac{1}{y_0y_{-1}y_{-2}}&0&0&0&0&0&\hdots
\\
\hdots&0&\dfrac{1}{y_0y_{-1}}&0&0&0&0&\hdots
\\
\hdots&0&0&\dfrac{1}{y_0}&0&0&0&\hdots
\\
\hdots&0&0&0&1&0&0&\hdots
\\
\hdots&0&0&0&0&y_1&0&\hdots
\\
\hdots&0&0&0&0&0&y_1y_2&\hdots
\\
\hdots&\hdots&\hdots&\hdots&\hdots&\hdots&\hdots&\ddots
\end{pmatrix}.
$$
\end{remark}

\begin{example}
\label{ex2.3}
In order to get the $\tau$-function~$H^\circ$ for simple Hurwitz numbers,
it suffices to set $y_c=e^{uc}$, $c=\dots,-2,-1,0,1,2,\dots$\,.
We prove this fact below in Sec.~\ref{ssec2.2}.
\end{example}

\subsection{Character formula for Hurwitz numbers}
\label{ssec2.2}
As we have mentioned already, Hurwiitz numbers are defined in terms of
combinatorics of symmetric groups. Consider this correspondence
in more detail. Let $S_n$~be the group of permutations of~$n$ elements,
let $\mathbb{C} S_n$ be its group algebra, and let $Z\mathbb{C} S_n$
be the center of the group algebra. For an arbitrary partition~$\mu$
of~$n$ denote by $C_\mu\in Z\mathbb{C} S_n$ the sum of the permutations
of cyclic type~$\mu$. All the permutations of the same cyclic type~$\mu$
are conjugate to one another and form a conjugacy class; let~$|C_\mu|$
denote the number of elements in this class. The elements $C_\mu/|C_\mu|$
form an additive basis in the vector space~$Z\mathbb{C} S_n$.
By definition, the Hurwitz number $h^\circ_{m;\mu}$ is the coefficient
of $C_\mu/|C_\mu|$ in the decomposition of $(C_2)^m/n!$ with respect to this basis;
here $C_2=C_{1^{n-2}2^1}$ is the sum of all transpositions.

In order to compute the element~$(C_2)^m$ and  for more general computations
in the algebra~$Z\mathbb{C} S_n$ it is more convenient, however, to use another
basis consisting of elements~$F_\lambda$, which form a system of idempotents
$F_\lambda F_\mu=\delta_{\lambda,\mu}F_{\lambda}$. By definition,
$F_\lambda\in Z\mathbb{C} S_n$ is the element that acts identically in the
irreducible representation~$\lambda$ and trivially in all the other irreducible
representations of~$S_n$. Recall that any central element acts as a scalar
operator in each irreducible representation. It follows that multiplication
by any element $a\in Z\mathbb{C} S_n$ is diagonal in the basis~$F_\lambda$:
$$
a=\sum_{\lambda\vdash n}f_a(\lambda)F_\lambda,\qquad
a\,F_\lambda=f_a(\lambda)F_\lambda.
$$
The number $f_a(\lambda)$ can be also defined as the constant multiplication by which
the element $a\in Z\mathbb{C} S_n$ acts in the irreducible representation~$\lambda$.
The theory of characters implies the following useful identity.

\begin{proposition}
\label{pr2.4}
For any element
$$
a=\sum_{\mu\vdash n}a_\mu\frac{C_\mu}{|C_\mu|}\in Z\mathbb{C} S_n,\qquad
a_\mu\in\mathbb{C},
$$
we have the following equation (in the space of polynomials in the variables
$p_1,p_2,\dots$):
$$
\sum_{\mu\vdash n}a_\mu p_{\mu_1}p_{\mu_2}\cdots=
\sum_{\mu\vdash n}f_a(\mu)\dim_\mu s_\mu(p),
$$
where $f_a(\mu)$ is the eigenvalue of the multiplication by~$a$ operator
corresponding to the basic eigenvector~$F_\mu$.
\end{proposition}

The latter identity is often used as a definition of Schur polynomials.
Note that in the standard interpretation of Schur polynomials as
symmetric functions in a set of variables $x_1,\dots,x_N$,
the variables~$p_k$ correspond to sums of the $k$~th powers of these variables
$p_k=x_1^k+\dots+x_N^k$.

In particular, applying this statement to the definition of Hurwitz numbers
we obtain the equation
$$
H^\circ(u;p_1,p_2,\dots)=\sum_{\lambda}e^{f_2(\lambda)u}
\frac{\dim_\mu}{|\mu|!}s_\lambda(p),
$$
where $f_2(\lambda)=f_{1^{n-2}2^1}(\lambda)$
is the corresponding eigenvalue of the operator of multiplication by the sum of transpositions.

In order to identify the function~$H^\circ$ completely, we have to compute the eigenvalues~$f_2(\lambda)$.
To this end, a formalism based on Jucys--Murphy elements is convenient.

\begin{definition}
\label{def4}
The \textit{Jucys--Murphy elements} of the group algebra~$\mathbb{C} S_n$
are the following sums of transpositions:
$$
X_1=0,\quad
X_2=(1,2),\quad
X_3=(1,3)+(2,3),\quad
X_4=(1,4)+(2,4)+(3,4),\quad\ldots,
$$
i.e.,
$$
X_k=\sum_{i=1}^{k-1}(i,k),\qquad
k=1,2,\dots,n.
$$
\end{definition}

These elements of the group algebra are not central.
Nevertheless, the following statements are true.

\begin{proposition}[\cite{14},~\cite{20},~\cite{29}]
\label{pr2.5}
{\rm1}. The elements $X_1,X_2,\dots,X_n$ commute pairwise, and hence
generate a commutative subalgebra in~$\mathbb{C} S_n$.

{\rm2}. Symmetric polynomials in~$X_k$ belong to the center
$Z\mathbb{C} S_n$, and each central element can be represented
as a symmetric polynomial in the Jucys--Murphy elements.

{\rm3}. If an element $a\in Z\mathbb{C} S_n$ is represented as a
symmetric polynomial in the Jucys--Murphy elements,
$a=P_a(X_1,\dots,X_n)$, then the eigenvalue~$f_a(\lambda)$
of the action of~$a$ in a representation~$\lambda$ is equal
to the value of the symmetric function~$P$ on the content~$c(\lambda)$
of the Young diagram~$\lambda$:
\begin{equation}
\label{eq6}
P_a(X_1,\dots,X_n)F_\lambda=P(c(\lambda))F_\lambda.
\end{equation}
\end{proposition}

Thus, for the element $C_2=\displaystyle\sum_{k=1}^nX_k$,
the eigenvalue corresponding to the Young diagram~$\lambda$ is
$$
f_2(\lambda)=\sum_{w\in\lambda}c(w)=
\frac{1}{2}\sum_{i=1}^{\ell(\lambda)}\biggl(\!\biggl(\lambda_i-i+
\frac{1}{2}\biggr)^2-\biggl(-i+\frac{1}{2}\biggr)^2\biggr).
$$
Hence, the generating function~$H^\circ$ is the element of the Orlov--Shcherbin
family with~$y_c=e^{c u}$, which proves the statement in Example~\ref{ex2.3}.

\begin{example}[{\rm(generalized Hurwitz numbers)}]
\label{ex2.6}
Simple Hurwitz numbers enumerate decompositions of a permutation of a given
cyclic type into a product of a given number of transpositions.
Consider a more general problem, where transpositions are replaced with
permutations having a given number of cycles. These numbers were introduced in~\cite{8}.

By the \textit{degeneracy of a partition}
$\lambda=(\lambda_1,\dots,\lambda_\ell)$ we mean the number
$$
k(\lambda)=|\lambda|-\ell(\lambda)=\sum_{i=1}^\ell(\lambda_i-1).
$$
The \textit{degeneracy of a permutation} is the degeneracy of the tuple
of lengths of independent cycles in it. In other words, the degeneracy of a permutation
is the minimal number of transpositions in its representation as a product of
transpositions. Set
$$
a^\circ_{k_1,\dots,k_m;\mu}=\frac{1}{|\mu|!}\bigl|\{(\tau_1,\dots,\tau_m), \
\tau_i\in S_{|\mu|}\colon k(\tau_i)=k_i, \
\tau_1\circ\dots\circ\tau_m\in C_\mu(S_{|\mu|})\}\bigr|
$$
and define the numbers~$a_{k_1,\dots,k_m;\mu}$ in a similar way
with an additional requirement that the subgroup in~$S_{|\mu|}$
generated by the permutations $\tau_1,\dots,\tau_\mu$
acts on the set $\{1,\dots,|\mu|\}$ transitively. The generating functions
\begin{equation}
\label{eq7}
A_m^\circ(u_1,\dots,u_m;p_1,p_2,\dots)=\sum_{\mu,k_1,\dots,k_m}
a^\circ_{k_1,\dots,k_m;\mu}u_1^{k_1}\cdots
u_m^{k_m}p_{\mu_1}p_{\mu_2}\cdots
\end{equation}
and
\begin{equation}
\label{eq8}
A_m(u_1,\dots,u_m;p_1,p_2,\dots)=\sum_{\mu,k_1,\dots,k_m}
a_{k_1,\dots,k_m;\mu}u_1^{k_1}\cdots u_m^{k_m}p_{\mu_1}p_{\mu_2}\cdots
\end{equation}
are related by
$$
A_m^\circ=\exp\{A_m\},
$$
and the following statement is true.

\begin{theorem}[\cite{8}]
\label{th2.7}
The generating function $A_m^\circ(u_1,\dots,u_m;p_1,p_2,\dots)$
for the generalized Hurwitz numbers is an $m$-parametric family of
$\tau$-functions for the KP hierarchy, and its logarithm
$A_m(u_1,\dots,u_m;p_1,p_2,\dots)$ is an
$m$-parametric family of solutions to this hierarchy.
\end{theorem}

The generating function~$A_m^\circ$ corresponds to the following set of
values of the parameters in the Orlov--Shcherbin family:
$$
y_c=\prod_{i=1}^m(1+u_ic)\quad\text{for}\quad
c=\dots,-2,-1,0,1,2,\dots,
$$
\begin{equation}
\label{eq9}
A_m^\circ(u_1,\dots,u_m;p_1,p_2,\dots)=\sum_\mu\biggl(\,\prod_{w\in\mu}\,
\prod_{i=1}^m(1+u_i c(w))\biggr)\frac{\dim_\mu}{|\mu|!}s_\mu(p_1,p_2,\dots).
\end{equation}

Indeed, denote by $C^{(k)}\in Z\mathbb{C} S_n$
the sum of all permutations of degeneracy~$k$.
Then $a^\circ_{k_1,\dots,k_m;\mu}$ is the coefficient of~$C_\mu/|C_\mu|$
in the decomposition of $(1/n!)C^{(k_1)}\cdots C^{(k_m)}$ in the basis~$C_\mu$.
Hence, Proposition~\ref{pr2.4} reduces computation of the generating function
for these numbers to the computation of the eigenvalues of the operator
of multiplication by~$C^{(k)}$ in~$Z\mathbb{C} S_n$. In order to compute these
eigenvalues, note that~$C^{(k)}$ can be expressed as the $k$~th elementary
symmetric function in the Jucys--Murphy elements:
\begin{equation}
\label{eq10}
C^{(k)}=\sigma_k(X_1,\dots,X_n),\qquad
\sum_{k=0}^n C^{(k)}u^k=\prod_{i=1}^n(1+X_i u).
\end{equation}
This statement follows from the fact that each permutation~$\tau$
admits a unique representation as a product of transpositions
$\tau=(a_1,b_1)\cdots(a_k,b_k)$, $a_i<b_i$, satisfying the following
\textit{strong monotonicity condition} for the transpositions entering it:
$b_1<\dots<b_k$. (It suffices to note that only the largest element
moved by~$\tau$ can be taken for~$b_k$, and argue by induction.)
Moreover, the number~$k$ of transpositions in such a product is exactly
the degeneracy of~$\tau$. Applying Eq.~\eqref{eq6}
from Proposition~\ref{pr2.5} we obtain Eq.~\eqref{eq9}
for the generating function for generalized Hurwitz numbers.
\end{example}

\begin{remark}
\label{rem2.8}
The above formula for the eigenvalues of
$\displaystyle\sum_{k=0}^n C^{(k)}u^k\in Z\mathbb{C} S_n$
can also be written as an identity for specializations of Schur functions
having one of the following two equivalent forms:
\begin{gather*}
s_\mu(v,v,v,\dots)=\frac{\dim_\mu}{|\mu|!}\prod_{w\in\mu}(v+c(w)),
\\
s_\mu(p_1,p_2,\dots)\big|_{p_i=u^{-1}}=u^{|\mu|}s_\mu(u^{-1},u^{-1},\dots)=
\frac{\dim_\mu}{|\mu|!}\prod_{w\in\mu}(1+u\,c(w)).
\end{gather*}
\end{remark}

\begin{example}[{\rm (Bousquet-M\'elou--Schaeffer numbers)}]
\label{ex2.9}
Similarly to the generalized Hurwitz numbers, the Bousquet-M\'elou--Schaeffer
numbers enumerate decompositions of a permutation of given cyclic type into
a product of a given number of permutations of arbitrary cyclic type.
In contrast to the generalized Hurwitz numbers, however, the total degeneracy
of the permutations in question rather than the degeneracy of each of them
is prescribed. These numbers were introduced in~\cite{5}.

Let us set
$$
b^\circ_{m,k;\mu}=\frac{1}{|\mu|!}\biggl|\biggl\{(\tau_1,\dots,\tau_m), \
\tau_i\in S_{|\mu|}\colon\sum_{i=1}^m k(\tau_i)=k, \
\tau_1\circ\dots\circ\tau_m\in C_\mu(S_{|\mu|})\biggr\}\biggr|
$$
and define the numbers~$b_{m,k;\mu}$ in a similar way, but with an additional
requirement that the subgroup in~$S_{|\mu|}$ generated by the permutations
$\tau_1,\dots,\tau_\mu$ must act on the set $\{1,\dots,|\mu|\}$ transitively.
The generating functions
\begin{align}
\label{eq11}
B_m^\circ(u;p_1,p_2,\dots)&=\sum_{k,\mu}
b^\circ_{m,k;\mu}u^kp_{\mu_1}p_{\mu_2}\cdots,
\\
\label{eq12}
B_m(u;p_1,p_2,\dots)&=\sum_{k,\mu\ne\varnothing}
b_{m,k;\mu}u^kp_{\mu_1}p_{\mu_2}\cdots
\end{align}
a related to one another by the equation
$$
B_m^\circ=\exp\{B_m\}.
$$
Obviously, they can be obtained from the generating functions for
generalized Hurwitz numbers by setting $u_1=\dots=u_m=u$:
\begin{align*}
B_m^\circ(u;p_1,p_2,\dots)&=A_m^\circ(u,\dots,u;p_1,p_2,\dots),
\\
B_m(u;p_1,p_2,\dots)&=A_m(u,\dots,u;p_1,p_2,\dots).
\end{align*}
This remark implies the following statement.

\begin{theorem}[\cite{8}]
\label{th2.10}
The generating function~$B_m^\circ(u;p_1,p_2,\dots)$
for Bousquet-M\'elou--Schaeffer numbers is a one-parameter family of
$\tau$-functions for the KP hierarchy, and its logarithm
$B_m(u;p_1,p_2,\dots)$ is a one-parameter family of solutions to this hierarchy.
\end{theorem}

Hence the generating function~$B_m^\circ$ corresponds to the following
set of values of the parameters
$$
y_c=(1+u c)^m\quad\text{for}\quad
c=\dots,-2,-1,0,1,2,\dots\quad
$$
in the Orlov--Shcherbin family:
$$
B_m^\circ(u;p_1,p_2,\dots)=\sum_\mu\,\prod_{w\in\mu}(1+u c(w))^m
\frac{\dim_\mu}{|\mu|!}s_\mu(p_1,p_2,\dots).
$$
\end{example}

\begin{remark}
\label{rem2.11}
If one does not impose the conditions on the degeneracy of the permutations
in the definition of Bousquet-M\'elou--Schaeffer numbers, that is, if we set
$u=1$ in the generating function for these numbers, then the answer becomes
extremely simple:
$$
B_m^\circ(1;p_1,p_2,\dots)=\sum_{k=0}^\infty (k!)^{m-1}s_k(p).
$$
Indeed, for $u=1$, the product $\displaystyle\prod_{w\in\mu}(1+u c(w))$
is nonzero only in the case when the Young diagram~$\mu$
consists of a single row, and it equals~$|\mu|!$ in this case.
\end{remark}

\begin{example}[{\rm (monotonic simple Hurwitz numbers)}]
\label{ex2.12}
The \textit{monotonic simple Hurwitz numbers} are defined as follows~\cite{7}.
Denote by~$\vec h^\circ_{m;\mu}$
(respectively,~$\vec h_{m;\mu}$) the result of dividing by~$|\mu|!$
the number of $m$-tuples of transpositions
$\tau_1,\dots,\tau_m$ in~$S_{|\mu|}$ possessing the following properties:

-- the cyclic type of the permutation $\tau_1\circ\dots\circ\tau_m$ is~$\mu$
(and, respectively, the subgroup generated by these transpositions acts
on the set $\{1,\dots,n\}$ transitively);

-- after writing each transposition~$\tau_i$ as the pair of elements
$(a_i,b_i)$ it permutes, so that $a_i<b_i$,
then $b_1\leqslant b_2\leqslant\dots\leqslant b_m$
(the \textit{monotonicity property}).

Similarly to the case of ordinary simple Hurwitz numbers, introduce two generating functions
\begin{align}
\label{eq13}
\vec H^\circ(u;p_1,p_2,\dots)&=\sum_{m=0}^\infty\,\sum_\mu
\vec h^\circ_{m;\mu}p_{\mu_1}p_{\mu_2}\cdots u^m,
\\
\label{eq14}
\vec H(u;p_1,p_2,\dots)&=\sum_{m=1}^\infty\,\sum_{\mu\ne\varnothing}
\vec h_{m;\mu}p_{\mu_1}p_{\mu_2}\cdots u^m.
\end{align}

Then the following statement is true.

\begin{theorem}[\cite{7},~\cite{13}]
\label{th2.13}
The generating function $\vec H^\circ(u;p_1,p_2,\dots)$
for simple monotonic Hurwitz numbers is a one-parameter family
of $\tau$-functions for the KP hierarchy, and its logarithm
$\vec H(u;p_1,p_2,\dots)$ is a one-parameter family of solutions
to this hierarchy.
\end{theorem}

The corresponding values of the parameters~$y_c$ in the family are
$$
y_c=\frac{1}{1-u c}\,,\qquad
c=\dots,-2,-1,0,1,2,\dots\,.
$$
Indeed, consider the sum of all possible monotonic decompositions of length~$m$
in the group~$S_n$:
$$
\vec C^{(m)}=\sum_{\substack{(i_1,\dots,i_m),(j_1,\dots,j_m)\\
i_1<j_1,\dots,i_m<j_m\\j_1\leqslant \dots\leqslant j_m}}
(i_1,j_1)\cdots(i_m,j_m).
$$
On one hand, the monotonic simple Hurwitz numbers are the coefficients
in the decomposition of $(1/n!)\vec C^{(m)}$ in the basis~$C_\mu/|C_\mu|$,
hence the generating function for these numbers is determined by the eigenvalues
of the operator of multiplication by~$\vec C^{(m)}$.
On the other hand, it is easy to see that~$\vec C^{(m)}$
can be represented as the $m$~th complete symmetric function in the Jucys--Murphy elements:
$$
\vec C^{(m)}=h_m(X_1,\dots,X_n),\qquad
\sum_{m=0}^\infty \vec C^{(m)}u^m=\prod_{k=1}^n\frac{1}{1-u X_k}\,,
$$
which implies the specialization of the parameters for the generating
function for monotonic simple Hurwitz numbers given above.

It is worth noting that the given specialization of the parameters in the
Orlov--Shcherbin family implies that the monotonic Hurwitz numbers correspond,
up to a sign, to the formal substitution ${m=-1}$ in the Bousquet-M\'elou--Schaeffer numbers:
$$
\vec h_{m,\mu}=(-1)^m b_{-1,m;\mu},\qquad
\vec H(u,p_1,p_2,\dots)=B_{-1}(-u,p_1,p_2,\dots).
$$
\end{example}

\begin{remark}[{\rm (Kerov--Olshanski algebra)}]
\label{rem2.14}
In the analysis of the character formula above, the number of permuted elements~$n$
for the group~$S_n$ is assumed to be given, in spite of the fact that
numbers corresponding to different values of~$n$ are collected into generating
series. S.~V.~Kerov and G.~I.~Olshanski~\cite{18}
suggested a natural way to unify the algebras $Z\mathbb{C} S_n$
for different values of~$n$ into a single algebra thus establishing a regular
dependence of the combinatorial quantities under consideration on~$n$.

Consider a countable set, say $\mathbb{N}$. And consider the set of all pairs
of the form (a finite subset in~$\mathbb{N}$, a permutation of this subset).
There is a natural semigroup operation on this set: we take the union of the two subsets,
and the product of the two permutations on this union.
The \textit{Kerov--Olshanski algebra~$\mathscr{W}$} is the center of the group
algebra of this semigroup. For any partition~$\mu$, denote by~$C_\mu$
the element of~$\mathscr{W}$ equal to the sum of all pairs of the form
(a $|\mu|$-element subset in~$\mathbb{N}$, a permutation of cyclic type~$\mu$
of this subset). Then the elements~$C_\mu$ for all partitions~$\mu$
form an additive basis in~$\mathscr{W}$. Here are some examples of multiplication in this algebra
(which is commutative and associative):
\begin{align*}
C_{1^1}^2&=C_{1^1}+2 C_{1^2},
\\
C_{1^1}C_{2^1}&=2C_{2^1}+C_{1^12^1},
\\
C_{2^1}^2&=3C_{3^1}+C_{1^2}+2C_{2^2},
\\
C_{2^1}C_{1^12^1}&=3C_{3^1}+4C_{2^2}+3C_{1^13^1}3C_{1^3}+2C_{1^12^2}.
\end{align*}

For any positive integer~$n$, there is a natural homomorphism
$\phi_n\colon\mathscr{W}\to Z\mathbb{C} S_n$,
which ``forgets the support of the permutation'':
$$
\phi_n(C_\mu)=\begin{pmatrix}
n-|\mu|+\varepsilon(\mu)
\\
\varepsilon(\mu)
\end{pmatrix}C_{1^{n-|\mu|}\mu},
$$
where $\varepsilon(\mu)$ is the number of ``ones'' in the partition~$\mu$. In particular,
$$
\phi_n(C_\mu)=\begin{cases}
C_\mu, &\text{if} \ |\mu|=n,
\\
0, &\text{if} \ |\mu|>n.
\end{cases}
$$
These homomorphisms allow one to associate to each element
$a\in\mathscr{W}$ a function on the set of Young diagrams:
$\lambda\mapsto f_a(\lambda)$, where $f_a(\lambda)$
is the eigenvalue of $\phi_{|\lambda|}(a)\in Z\mathbb{C} S_{|\lambda|}$,
which corresponds to the partition~$\lambda$. The \textit{Kerov--Olshanski theorem}
means (in the modern language) that for each $a\in\mathscr{W}$
the function~$f_a$ on the set of Young diagrams is a symmetric polynomial
in the contents of the diagrams; moreover, this correspondence establishes
an isomorphism between the algebra~$\mathscr{W}$ and the algebra of symmetric functions.
Another version of the Kerov--Olshanski theorem interprets~$f_a(\lambda)$
as the value of a certain symmetric polynomial in the shifted parts
$\lambda_1-1,\lambda_2-2,\lambda_3-3,\dots$ of the partition
$\lambda=(\lambda_1,\lambda_2,\dots)$,
$\lambda_1\geqslant\lambda_2\geqslant\cdots$\,.
That is why the algebra~$\mathscr{W}$ often carries the name of the
\textit{algebra of shifted symmetric functions} (see~\cite{23}).

In the context of Hurwitz numbers, Proposition~\ref{pr2.5}
states that the value of any symmetric polynomial (in an unspecified number of elements)
on Jucys--Murphy elements defines an element of the Kerov--Olshanski algebra.
Regularity of the dependence of these numbers on~$n=|\mu|$ then follows.
\end{remark}

\subsection{Genus expansion}
\label{ssec2.3}
Hurwitz numbers are closely related to enumeration of ramified coverings of the projective
line. Namely, the Hurwitz numbers~$h^\circ_{m,\mu}$
enumerate ramified coverings of the projective line~$\mathbb{C}P^1$
with prescribed ramification points in it such that
the ramification over one of them has a given cyclic type~$\mu$,
while it is simple over the other ramification points.
Ramified coverings are enumerated with the weight inverse to the order of
the automorphism group of the covering (note, however, that the latter
group is usually trivial). The genus of the covering surface is uniquely determined
from the ramification data by the Riemann--Hurwitz formula; however, the
covering surface can be disconnected, and it is better to talk about its
Euler characteristic rather than its genus. In its own turn, the logarithm~$H(u;p_1,p_2,\dots)$
of the $\tau$-function~$H^\circ$ enumerates connected ramified coverings.
In the case of the connected ramified covering of the projective line
having ramification type~$\mu$ over one point and simple ramifications over~$m$
other points, the \textit{Riemann--Hurwitz formula for the genus~$g$ of the covering surface} acquires the form
$$
2g-2=m-|\mu|-\ell(\mu).
$$
It gives a hint how to extract the genus~$g$ of the connected covering surface
in the generating function explicitly: each additional point of simple ramification
must add~$1$ to the degree of the parameter, which is in charge of the genus,
and a part of length~$i$ in the partition~$\mu$ must diminish this degree by~$i+1$.
As a result, we must make the substitution
$p_i\mapsto \hbar^{-i-1} p_i$, which results in the deformation of the Schur polynomials:
$$
s_\lambda^\hbar(p_1,p_2,\dots)=s_\lambda\biggl(\frac{p_1}{\hbar^2}\,,
\frac{p_2}{\hbar^3}\,,\frac{p_3}{\hbar^4}\,,\dots\biggr),
$$
and set $y_c=e^{u\hbar c}$ in the parameters of the Orlov--Shcherbin family.
Then the generating function for connected simple Hurwitz numbers
acquires the following genus expansion:
\begin{align*}
H^\hbar(u;p_1,p_2,\dots)&=\hbar^2\log\biggl(\,\sum_\mu e^{u\hbar c(\mu)}
\frac{\dim_\mu}{|\mu|!}s^\hbar_\mu(p_1,p_2,\dots)\biggr)
\\
&=H_0(u;p_1,p_2,\dots)+\hbar^2H_1(u;p_1,p_2,\dots)
\\
&\qquad+\hbar^4H_2(u;p_1,p_2,\dots)+\cdots,
\end{align*}
where the function~$H_g$, $g=0,1,2,\dots$, enumerates ramified coverings of genus~$g$
In its own turn, the KP hierarchy also will be deformed. In particular, the first KP
equation~\eqref{eq1} for the deformed function becomes
$$
\frac{\partial^2 H^\hbar}{\partial p_2^2}=
\frac{\partial^2 H^\hbar}{\partial p_1\,\partial p_3}-
\frac{1}{2}\biggl(\frac{\partial^2 H^\hbar}{\partial
p_1^2}\biggr)^2-\frac{\hbar^2}{12}\,\frac{\partial^4 H^\hbar}{\partial p_1^4}
$$
(we make use of the identity
$H^\hbar(u;p_1,p_2,\dots)=\hbar^2H(\hbar u;p_1/\hbar^2,p_2/\hbar^3,\dots)$).
This means, for example, that the generating function~$H_0$
for connected Hurwitz numbers of genus~$0$ is a solution to the
\textit{dispersionless KP hierarchy}:
$$
\frac{\partial^2 H_0}{\partial p_2^2}=
\frac{\partial^2 H_0}{\partial p_1\,\partial p_3}-
\frac{1}{2}\biggl(\frac{\partial^2 H_0}{\partial p_1^2}\biggr)^2.
$$

The next term in the series expansion in the parameter~$\hbar$
satisfies the equation
$$
\frac{\partial^2 H_1}{\partial p_2^2}=
\frac{\partial^2 H_1}{\partial p_1\,\partial p_3}-
\frac{\partial^2 H_0}{\partial p_1^2}\,\frac{\partial^2 H_1}{\partial p_1^2}
-\frac{1}{12}\,\frac{\partial^4 H_0}{\partial p_1^4}\,,
$$
which, in contrast to the original KP equation, is linear with respect to the
unknown function~$H_1$. Therefore, if we know the zero genus expansion term~$H_0$,
then we can find the next term~$H_1$ by means of standard methods of the theory
of linear differential equations. The same is true for all higher genera.

In the case of simple Hurwitz numbers, the zero genus term~$H_0$
is known since Hurwitz. Its coefficients are
$$
h_{m;\mu}=|\mu|^{\ell(\mu)-3}\frac{m!}{|\operatorname{Aut}(\mu)|}
\prod_{i=1}^{\ell(\mu)}\frac{\mu_i^{\mu_i}}{\mu_i!}\,,\qquad
m=|\mu|+\ell(\mu)-2.
$$

Similarly to the ordinary Hurwitz numbers, all their modifications also
can be interpreted as numbers of ramified coverings of the sphere satisfying
prescribed restrictions on the ramification types at given ramification points.
The transitivity condition for the subgroup generated by the permutations action in each case
is nothing but the connectivity condition for the covering surface.
The genus of this surface is determined from the Riemann--Hurwitz formula
and can be extracted explicitly by an appropriate renormalization of the parameters.
Moreover, we can introduce the notion of genus expansion for an arbitrary
solution of the KP hierarchy belonging to the Orlov--Shcherbin family.

Let
\begin{equation}
\label{eq12n}
\varphi(c)=d_0+d_1c+d_2c^2+d_3c^3+\dotsb
\end{equation}
be an arbitrary power series.

\begin{definition}
\label{def5n}
The \textit{genus expansion} of a $\tau$-function obtained from a function in the Orlov--Shcherbin family
$$
\sum_\mu y_\mu\frac{\dim_\mu}{|\mu|!}s_\mu(p_1,p_2,\dotsb)
$$
by a substitution of the form
$y_c=\varphi(c)$, $c=\dotsc,-2,-1,0,1,2$, is the function
$$
\sum_\mu \biggl(\prod_{w\in\mu} \varphi(\hbar c(w))\biggr)\frac{\dim_\mu}{|\mu|!}s^\hbar_\mu(p_1,p_2,\dotsc).
$$
\end{definition}

Let us prove that this notion is well-defined.

\begin{theorem}
\label{th2.15n}
For any power series~\eqref{eq12n} the power series expansion in the parameter~$\hbar$
of the function
$$
\Phi^\hbar(p_1,p_2,\dots)=\hbar^2\log \sum_\mu \biggl(\prod_{w\in\mu} \varphi(\hbar c(w))\biggr)\frac{\dim_\mu}{|\mu|!}s^\hbar_\mu(p_1,p_2,\dotsc)
$$
contains only even nonnegative powers of~$\hbar$.
\end{theorem}

Hence, each monomial in the power expansion indeed has its own genus.

It seems that this theorem must have an independent proof, not related to
ramified coverings of the sphere. Alternatively, it can be proved using a result
due to Guay-Paquet and Harnad~\cite{11}, which interprets the function
$\Phi^\hbar$ as the generating function for the numbers of ramified coverings
of the sphere by connected surfaces counted with certain weights. Consider
a (nonstrictly) monotonic decomposition of an $n$-element permutation into
a product of transpositions
$$
(i_1,j_1)(i_2,j_2)\dotsb(i_m,j_m),\qquad i_k<j_k, \quad j_1\le j_2\le\dotsb\le j_m.
$$
The \textit{weight} of such a permutation is a monomial $d_{k_1}\dotsb d_{k_n}$
in the formal variables $d_0,d_1,\dotsc$, where $k_s$
is the number of those transpositions in the decomposition whose second (the greater one) index is~$s$.
Respectively, the total number of transpositions in the decomposition is
$m=k_1+\dotsb+k_n$. Then the $\tau$-function in Definition~\ref{def5n}
is the generating function for the numbers of nonstrictly monotonic decompositions
counted with the weights, while the function $\Phi^\hbar$ of the theorem
is the generating function for similar decompositions satisfying the transitivity condition.
The power of the variable~$\hbar$ in each monomial in~$\Phi^\hbar$  is twice
the genus of the covering surface given by the decomposition.
Indeed, consider the product
$$\varphi(X_1)\varphi(X_2)\dotsb \varphi(X_n)$$
(as an element of the ring $Z\mathbb C S_n[[d_0,d_1,\dots]]$), where $X_1,\dots,X_n$ are the Jucys--Murphy element.
Then, as one can easily see, this product is the sum of all nonstrictly monotonic decompositions
into a product of transpositions in the group~$S_n$ (for a given~$n$),
taken with coefficients equal to the weights of the corresponding decompositions.
This remark implies the Guay-Paquet--Harnad theorem in the way similar to the proofs of
Theorems~\ref{th2.7},~\ref{th2.10} and~\ref{th2.13}
above.

\begin{example}
\label{ex16}
The generating function~$B_m$ for connected Bousquet-M\'elou--Schaeffer numbers
corresponds to the set of parameters $y_c=(1+u c)^m$
in the Orlov--Shcherbin family. In order to obtain its genus expansion~$B_m^\hbar$
it suffices to replace these parameters by $y_c=(1+u\hbar c)^m$:
$$
B_m^\hbar=\hbar^2\log\sum_\mu\biggl(\prod_{w\in\mu}(1+u\hbar c(w))^m\biggr)
\frac{\dim_\mu}{|\mu|!}s^\hbar_\mu(p_1,p_2,\dots).
$$
The coefficients of this function enumerate connected coverings of the projective line
by a surface of given genus with prescribed ramification points, over one of which
the ramification has the cyclic type~$\mu$, while all the other ramification points
are among the other~$m$ points, and the ramification type over each of them can be arbitrary.
The following elegant explicit formula has been obtained for the number of genus zero
coverings in~\cite{5}:
$$
b_{m,k,\mu}=\frac{m(m|\mu|-k+1)^{\overline{\ell(\mu)-3}}}
{|\operatorname{Aut}(\mu)|}\,\prod_{i=1}^{\ell(\mu)}
\begin{pmatrix} m\mu_i-1 \\ \mu_i \end{pmatrix},\qquad
k=|\mu|+\ell(\mu)-2,
$$
where
$$
(d+1)^{\overline{r}}=(d+1)(d+2)\cdots(d+r)
$$
denotes the ascending product of~$r$ factors, and for $r<0$ we set
$$
(d+1)^{\overline{r}}=\frac{1}{(d+r+1)^{\overline{-r}}}\,.
$$
\end{example}

\begin{example}
\label{ex2.16}
For monotonic Hurwitz numbers, the genus expansion has the form
$$
\vec H^\hbar=\hbar^2\log\sum_{\mu}\prod_{w\in\mu}\,\frac{1}{1-u\hbar c(w)}\,
\frac{\dim_\mu}{|\mu|!}s^\hbar_\mu(p_1,p_2,\dots).
$$
The formula for the monotonic Hurwitz numbers of genus~$0$ has been obtained in~\cite{6}.
Up to the sign~$(-1)^m$ it can be obtained by the formal substitution $m=-1$,
$k=m$ in the above formula for the genus~$0$ Bousquet-M\'elou--Schaeffer numbers:
$$
\vec h_{m,\mu}=\frac{(2|\mu|+1)^{\overline{\ell(\mu)-3}}}
{|\operatorname{Aut}(\mu)|}\prod_{i=1}^{\ell(\mu)}
\begin{pmatrix} 2\mu_i \\ \mu_i\end{pmatrix},\qquad
m=|\mu|+\ell(\mu)-2.
$$
\end{example}

\subsection{Toda lattice hierarchy}
\label{ssec2.4}
Theorem~\ref{th2.1} is a special case of the following more general result.
Consider two sequences of independent variables $p=(p_1,p_2,\dots)$ and $q=(q_1,q_2,\dots)$,
and an additional set of formal parameters
$(\dots,y_{-1},y_0,y_1,y_2,\dots)$. For each $n\in\mathbb{Z}$, set
$$
\tau_n(p_1,p_2,\dots;q_1,q_2,\dots)=r_0(n)\sum_\mu r_\mu(n)
s_\mu(p_1,p_2,\dots)s_\mu(q_1,q_2,\dots),
$$
where summation is carried over all partitions~$\mu$ and
$$
r_\mu(n)=\prod_{w\in\mu}y_{c(w)+n},\qquad
r_0(n)=\begin{cases}
\displaystyle\prod_{j=1}^{n-1} y_j^{n-j},&n\geqslant0,
\\
\displaystyle\prod_{j=n+1}^0 y_j^{j-n},&n<0.
\end{cases}
$$

\begin{theorem}
\label{th2.17}
Each function~$\tau_n$ is a $\tau$-function of the KP hierarchy both
in the set of variables~$p$ and in the set of variables~$q$.
\end{theorem}

In fact, a stronger statement is valid, namely, the set of functions~$\tau_n$, $n\in\mathbb{Z}$,
forms a set of $\tau$-functions of the Toda lattice hierarchy.
The Toda lattice hierarchy~\cite{28} contains the KP hierarchy equations
for $\log \tau_n$ with respect to each of the two groups of variables
and for each~$n$, as well as an additional set of equations including derivatives
with respect to variables from the two groups of functions with three subsequent indices.
The simplest such equation is
$$
\frac{\partial^2\log\tau_n}{\partial p_1\,\partial q_1}=
\frac{\tau_{n-1}\tau_{n+1}}{\tau_n^2}\,,
\qquad n\in\mathbb{Z}.
$$

Similarly to the KP hierarchy, it is more geometrically transparent to
describe the space of solutions to the Toda lattice hierarchy
rather than the equations themselves. The solutions are parameterized
by infinite dimensional matrices that are elements~$M$ of the group
$\operatorname{GL}(V)$ (see~\cite[Ch.~4]{15}). For the sake of simplicity,
we restrict ourselves to the case where~$M$ is an upper triangular matrix
$$
M\colon z^k\mapsto\sum_{i\geqslant k}M_{i,k}z^i,\qquad
M_{k,k}\ne 0.
$$
The solution corresponding to $M=\operatorname{Id}$ is
$$
\tau_n=\exp\biggl\{\,\sum_k\frac{p_kq_k}{k}\biggr\}=\sum_\mu s_\mu(p)s_\mu(q)
$$
for all~$n$. Under the boson--fermion correspondence (in $p$-variables)
such a function~$\tau_n$ is transformed into a decomposable wedge product of the form
$\beta_1\wedge\beta_2\wedge\cdots$, where
$$
\beta_k=\exp\biggl\{\,\sum_i\frac{q_iz^i}{i}\biggr\}z^{n-k}=
\sum_{j=0}^\infty s_j(q)z^{n-k+j}.
$$
For $n\ne0$, this wedge product does not belong to the space~$\Lambda^{\infty/2}V$;
instead, it is an element of another space, denoted by~$\Lambda^{\infty/2+n}V$ and
named the \textit{space of fermions of charge~$n$}. Similarly to the case of zero charge,
the elements in the standard basis of $\Lambda^{\infty/2+n}V$ are indexed by partitions.
They have the form
$$
v_\mu=z^{m_1}\wedge z^{m_2}\wedge z^{m_3}\wedge\cdots,\qquad
m_i=\mu_i-i+n.
$$

For a general element $M\in \operatorname{GL}(V)$, the fermion
associated to the $\tau$-function~$\tau_n$ undertakes the following modification:
\begin{equation}
\label{eq15}
\tau_n\leftrightarrow r_0(n)\frac{M(\beta_1)}{M_{n-1,n-1}}\wedge
\frac{M(\beta_2)}{M_{n-2,n-2}}\wedge\cdots,
\end{equation}
where
$$
r_0(n)=\begin{cases}
\displaystyle\prod_{j=0}^{n-1} M_{j,j},&n\geqslant0,
\\
\displaystyle\prod_{j=n}^{-1} M_{j,j}^{-1},&n<0.
\end{cases}
$$

The normalizing factor~$r_0(n)$ on the right-hand side of Eq.~\eqref{eq15}
is due to the fact that the vector spaces we consider are infinite dimensional.
This factor is independent of~$n$ and equals $\displaystyle\prod_{i=1}^\infty M_{-i,-i}$;
its parts are grouped so as to make the coefficient of the term
of the smallest degree in the $k$~th factor (that is, the coefficient of~$z^{n-k}$)
equal to~$1$ for~$k$ large enough. It is easy to see that the solution to the hierarchy
in Theorem~\ref{th2.17} corresponds to the case where $M$ is the diagonal matrix
in Remark~\ref{rem2.2}.

\begin{example}[{\rm (double Hirwitz numbers and Toda lattice hierarchy)}]
\label{ex2.18}
The statement that the generating function for double Hurwitz
numbers is a $\tau$-function for an integrable hierarchy proved
by A.~Okounkov in~\cite{22}
was the first one and served as a model for statements of this kind.

In order to define double Hurwitz numbers, consider
for a given number~$m$ all tuples of permutations
$(\alpha,\beta,\tau_1,\dots,\tau_m)$, where

-- the permutation~$\alpha$ belongs to the cyclic type specified
by a given partition $\mu\vdash|\mu|$ and, therefore, is an element
of the group~$S_{|\mu|}$;

-- the permutation~$\beta\in S_{|\mu|}$ belongs to the cyclic type
specified by a given partition $\nu\vdash|\nu|=|\mu|$;

-- the permutations $\tau_1,\dots,\tau_m$ are transpositions;

-- the permutation $\alpha\circ\beta\circ\tau_1\circ\dots\circ\tau_m$
is the identity.

The \textit{double Hurwitz number~$d^\circ_{m;\mu,\nu}$}
(respectively, the  \textit{connected double Hurwitz number~$d_{m;\mu,\nu}$})
is the number
\begin{align*}
d^\circ_{m;\mu,\nu}&=\frac{1}{|\mu|!}
|\{(\alpha,\beta,\tau_1,\dots,\tau_m)\colon
\alpha\in C_\mu(S_{|\mu|}), \ \beta\in C_\nu(S_{|\mu|}),
\\
&\qquad\tau_i\in C_2(S_{|\mu|}), \ \alpha\circ\beta\circ\tau_1\circ\dots\circ
\tau_m=\operatorname{id}\}|
\end{align*}
(in the definition of~$d_{m;\mu,\nu}$ we require in addition
that the action of the subgroup in~$S_{|\mu|}$
generated by the permutations $\alpha,\beta,\tau_1,\dots,\tau_m$ must be transitive).

Double Hurwitz numbers enumerate ramified coverings of the two-dimensional sphere
having ramification of given cyclic types~$\mu$ and~$\nu$
over two given critical values. Under an appropriate choice of a complex coordinate~$z$
in the covered sphere these points will acquire the coordinates $z=0$ and $z=\infty$.
Connected double Hurwitz numbers enumerate connected coverings.

Let us collect double and connected double Hurwitz numbers into two generating
functions, each depending on two infinite families of variables
$q_1,q_2,\dots$ and $p_1,p_2,\dots$, and additional parameters~$u$,~$v$:
\begin{align}
\label{eq16}
D^\circ(u,v;p_1,p_2,\dots;q_1,q_2,\dots)&=\sum_{m=0}^\infty\sum_{\mu,\nu}
d^\circ_{m;\mu,\nu}p_{\mu_1}p_{\mu_2}\cdots q_{\nu_1}q_{\nu_2}\cdots v^{|\mu|}
\frac{u^m}{m!}\,;
\\
\label{eq17}
D(u,v;p_1,p_2,\dots;q_1,q_2,\dots)&=\sum_{m=1}^\infty\sum_{\mu,\nu}
d_{m;\mu,\nu}p_{\mu_1}p_{\mu_2}\cdots q_{\nu_1}q_{\nu_2}\cdots v^{|\mu|}
\frac{u^m}{m!}\,.
\end{align}

\begin{theorem}[\cite{22}]
\label{th2.19}
The generating function~$D^\circ$ is a family of $\tau$-functions
for the Toda lattice hierarchy.
\end{theorem}
\end{example}

\begin{example}
\label{ex2.20}
Similarly to double Hurwitz numbers, the generalized double Hurwitz
numbers, the double Bousquet-M\'elou--Schaeffer numbers
and the monotonic double Hurwitz numbers can be defined.
Generating functions for these versions of double Hurwitz numbers
depend on two sets of time variables $(p_1,p_2,\dots)$
and $(q_1,q_2,\dots)$, and some additional variables. They can be obtained
from the function in Theorem~\ref{th2.17} by substituting the same
values of~$y_c$ as in the previous section in the study of simple numbers.
As a consequence, these generating functions for connected numbers
satisfy KP hierarchy with respect to each set of time variables
and Toda lattice hierarchy, and the corresponding generating functions for simple
numbers are results of the substitution $q_1=1$, $q_i=0$ ($i>1$)
in the generating functions for double numbers.
\end{example}

During the last years the set of examples of combinatorial solutions
to integrable hierarchies grows very quickly
(see, for example,~\cite{1}--\cite{3},~\cite{10}--\cite{13},
\cite{17},~\cite{21}). On the other hand, all attempts
to formulate and prove similar integrability properties
for triple and more general Hurwitz numbers, whose definition requires
specifying the complete ramification profile over more than two
points of the sphere, we know about were unsuccessful.

\begin{example}
\label{ex2.21}
Below, we will require the following family of solutions to
the Toda lattice hierarchy:
\begin{align*}
N(u;p_1,p_2,\dots;q_1,q_2,\dots)&=\log\sum_{|\mu|}\,\prod_{w\in\mu}
(u+c(w))s_\mu(p_1,p_2,\dots)s_\mu(q_1,q_2,\dots)
\\
&=\sum_{|\lambda|=|\mu|}\,\sum_{k=1}^{|\mu|}
n_{k;\mu,\lambda}u^k p_\lambda q_\mu.
\end{align*}
The second equation is the definition of the coefficients~$n_{k;\lambda,\mu}$.
We conclude that~ $|\lambda|!\,n_{k;\lambda,\mu}$ (where $|\lambda|=|\mu|$)
is the number of triples of permutations~$\sigma$,~$\alpha$,~$\varphi$
of~$|\lambda|$ elements satisfying the following conditions:

-- the product $\varphi\alpha\sigma$ is the identity permutation;

-- the permutations~$\alpha$ and~$\sigma$ are of cyclic types~$\mu$
and~$\lambda$, respectively;

-- the number of cycles in the permutation~$\sigma$ is~$k$;

-- the group generated by the permutations~$\varphi$,~$\alpha$,~$\sigma$,
acts transitively on the set of permuted elements.

Note that in spite of the obvious symmetry between the permutations~$\varphi$, $\alpha$,~$\sigma$,
the cyclic type of only two of them is taken into account in the construction of the generating function~$N$,
while only the number of cycles is extracted from the third permutation.
\end{example}

\subsection{Applications: enumeration of maps and triangulations}
\label{ssec2.5}
A map on a given surface is a graph embedded into the surface and cutting
it into two-dimensional discs. We will assume that the surface is compact and oriented.
Two maps are said to be equivalent if there is an orientation preserving
homeomorphism of one surface onto the second one taking the vertices and
the edges of the first graph to the vertices and the faces of the second one.
Gluing a surface from a set of polygons provides an alternative way of looking at a map:
the graph then is the result of gluing the polygons' boundaries.

For the purposes of the present paper it will be more convenient for us to
use the following combinatorial definition of a map.

\begin{definition}
\label{def5}
Let~$D$ be a finite set. A \textit{map with the set of half-edges}~$D$
\textit{on an oriented surface} is a triple of permutations~$\alpha$,~$\varphi$,~$\sigma$
of the set~$D$ possessing the following properties:

-- $\alpha$ is an involution without fixed points;

-- the product $\varphi\alpha\sigma$ is the identity permutation.

The group $G=\langle\alpha,\varphi,\sigma\rangle$ of permutations of the set~$D$
generated by the permutations~$\alpha$,~$\varphi$,~$\sigma$,
is called the \textit{cartographic group} of the map.
The map is said to be \textit{connected} if its cartographic group acts
on the set~$D$ transitively.
\end{definition}

In the case of a graph drawn on an oriented surface, $D$~is the set
of half-edges, or flags, of the graph, the permutation~$\alpha$
exchanges the ends of each edge, the permutation~$\varphi$
rotates the half-edges around each face, and~$\sigma$
rotates the half-edges around each vertex (see details in~\cite{33}).
Obviously,~$\alpha$ is an involution without fixed points, and it is easy
to verify that the product of these three permutations indeed is
the identity permutation. A map is connected iff the surface where it is drawn is connected.

The number of edges in the map is half the number of elements in the set~$D$,
that is, it is equal to the number of cycles in the permutation~$\alpha$.
The number of vertices in the map is the number of cycles in~$\sigma$,
the valencies of the vertices are the lengths of the cycles.
Similarly, the number of faces in the map is the number of cycles in
the permutation~$\varphi$, and the degrees of the faces are the lengths of the cycles.

The notion of \textit{hypermap} generalizes that of the map.
In the definition of hypermap we get rid of the restriction that~$\alpha$
is an involution without fixed points, thus restoring the symmetry between
the three permutations.

Now it is clear that enumeration of maps and hypermaps of various kinds
is reduced to enumeration of triples of permutations possessing certain
special properties, and hence the enumerative methods described above can be applied.
In particular, Example~\ref{ex2.21} shows that the generating function
whose coefficients enumerate hypermaps with given sets of valencies of
vertices and faces is a family of solutions to the Toda lattice hierarchy.
As a result, this function, as well as various its specializations,
including those that enumerate maps, are solutions to the KP hierarchy.

Thus, denote by~$R_\kappa^{(n,m)}$ the number of connected rooted maps with~$n$
edges, $m$~faces and the valencies of the verticies given by a partition~$\kappa$ of~$2n$.
Then the following theorem is true.

\begin{theorem}[\cite{8}]
\label{th2.22}
The generating function
$$
R(w,z;p_1,p_2,\dots)=\sum_{n,m\geqslant1}\,\sum_{\kappa\vdash 2n}
\frac{R_\kappa^{(n,m)}}{2n}p_\kappa w^mz^n
$$
(where, for a given partition $\kappa=(\kappa_1,\kappa_2,\kappa_3,\dots)$,
$p_\kappa$ denotes the monomial
$p_\kappa=p_{\kappa_1}p_{\kappa_2}p_{\kappa_3}\cdots$)
is a two-parameter family of solutions to KP equations.
\end{theorem}

Indeed, in contrast to the number $R_\kappa^{(n,m)}$ of rooted maps
the ratio ${R_\kappa^{(n,m)}}/{(2n)}$ is a weighted (i.e., counted with
automorphisms taken into account) number of ordinary (unrooted) maps,
that is, it is equal to the number $a_{m,n}(\kappa)$ in Example~\ref{ex2.6}.
Therefore, the function~$R$ of the theorem coincides with the function~$A_m$
of Theorem~\ref{th2.7} for the case $m=2$, $u_1=w$, $u_2=z$
(alternatively, the function~$R$ is a specialization of the function~$N$
in Example~\ref{ex2.21} for $u=w$, $q_i=z^i$, $i=1,2,\dotsc$).

For certain values of the parameters in the series~$R$ we obtain
the generating function for rooted cubic maps, that is, maps
valencies of all whose verticies are equal to~$3$. By duality,
this is equivalent to enumeration of rooted triangulations of surfaces
of arbitrary genus. The KP equations then provide recurrent relations
for the numbers of rooted triangulations.

Denote the number of rooted triangulations of a surface of genus~$g$
with~$2n$ faces by~$T(n,g)$. Then the recurrent relations acquire the
form below. Introduce notation
$$
S=\biggl\{(n,g)\in\mathbb{Z}\times\mathbb{Z}\colon n\geqslant-1,\
0\leqslant g\leqslant\frac{n+1}{2}\biggr\}.
$$

\begin{theorem}[\cite{8}]
\label{th2.23}
We have
$$
T(n,g)=\frac{1}{3n+2}t(n,g),
$$
where the numbers $t(n,g)$ are defined by the quadratic recurrent relation
$$
t(n,g)=\frac{4(3n+2)}{n+1}\biggl(n(3n-2)t(n-2,g-1)+
\sum t(i,h)t(j,k)\biggr)
$$
for $(n,g)\in S\setminus\{(-1,0),(0,0)\}$, where the sum is taken over elements
$(i,h)\in S$, $(j,k)\in S$ such that $i+j=n-2$ and
$h+k=g$, under the initial condition
$$
t(-1,0)=\frac{1}{2}\,,\qquad
t(n,g)=0\quad\text{for} \
(n,g)\notin S.
$$
\end{theorem}

This recurrent relation allowed E.~A.~Bender, Z.~Gao, and L.~B.~Richmond
to give a purely combinatorial computation of the constant factor in the
leading term in the asymptotics of the number of rooted triangulations
as the number of triangles tends to infinity. This coefficient was previously
computed using either matrix models of quantum gravity~\cite{30},
or intersection theory on moduli spaces~\cite{15n}.

\begin{theorem}[\cite{4}]
\label{th2.24}
The number of rooted triangulations of a genus~$g$ surface consisting of~$2n$
triangles has the asymptotics
$$
T(n,g) \sim \frac{3\,b_g}{\Gamma(5g/2-1/2)} \biggl(\frac38\biggr)^{(g-1)/2}n^{5(g-1)/2}(12\sqrt{3})^n\quad
\text{for } n\to\infty.
$$
The constant~$b_g$ here is defined by the initial condition
$b_0=-1$ and the quadratic recurrent equation
$$
b_{g+1}=\frac{25g^2-1}{24}b_{g}+\frac12\sum_{m=1}^{g}b_{g+1-m}b_{m} \quad \text{for } g\ge1.
$$
\end{theorem}

The recurrent equation is equivalent to the statement that the generating
function for the constants $b_g$ is a solution to the Painlev\'e~I equation:
$$U(y)=\sum_{g=0}^\infty\frac{b_g}{y^{5(g-1)/2+2}},\qquad \frac13 U''+U^2=y.$$
The first few values of the constant~$b_g$ are
\begin{gather*}
b_0=-1,\qquad b_1=\frac1{24},\qquad
b_2=\frac{49}{1152},\qquad
b_3=\frac{1225}{6912}, \qquad b_4=\frac{4412401}{2654208}.
\end{gather*}
This constant governs many other similar asymptotics.
In particular, it appears in~\cite{32} in the description of a
two-dimensional quantum gravity model based on simple Hurwitz numbers.
Namely, the following statement is true.

\begin{theorem}
\label{th2.25}
For a given value~$g$ of the genus, we have the equivalence
$$
\frac{h_{2n+2g-2;1^n}}{(2n+2g-2)!}\sim e^nn^{5(g-1)/2-1}\frac{b_g}{\Gamma\bigl(5g/2-1/2\bigr)2^{3g/2-1/2}}
$$
as $n\to\infty$.
\end{theorem}

The constant $b_g$ also can be explicitly expressed in terms of intersection numbers
on the moduli spaces of stable curves of genus~$g$ with~$n$ marked points
$\overline{{\mathcal M}}_{g;n}$
(see, e.g., \cite{15n},~\cite{22n})
$$
b_g=\frac{(5g-5)(5g-3)}{2^g(3g-3)!}\int_{\overline{{\mathcal M}}_{g;3g-3}}\psi_1^2\psi_2^2\dotsb\psi_{3g-3}^2.
$$

Mention also a recent result due to P.~G.~Zograf related to enumeration of
dessins d'enfants and having the same origin.

Denote by~$b_{m,n}(\kappa)$ the number of rooted Belyi graphs
corresponding to a partition~$\kappa$ with~$m$ vertices and~$n$ faces, and let
$$
B_d(u,v;p_1,p_2,\dots)=\sum b_{m,n}(\kappa) u^mv^np_\kappa,
$$
where the summation is carried over all~$m$,~$n$ and $\kappa\vdash d$ \ldots\,.
Set
$$
\mathscr{B}(s,u,v;p_1,p_2,\dots)=\sum_{d\geqslant0}B_d(u,v;p_1,p_2,\dots)s^d.
$$

\begin{theorem}[\cite{31}]
\label{th2.26}
The generating function $\mathscr{B}(s,u,v;p_1,p_2,\dots)$
is a three-parameter family of solutions to the KP hierarchy.
\end{theorem}


\end{document}